\documentclass[reqno,11pt]{article}

\usepackage{amsmath,latexsym,amssymb,theorem}

\def\sqr#1#2{{\vcenter{\hrule height.#2pt
        \hbox{\vrule width.#2pt height#1pt \kern#1pt
                \vrule width.#2pt}
        \hrule height.#2pt}}}

\newtheorem{Theorem}{\sc Theorem}[section]
\newtheorem{Lemma}[Theorem]{\sc Lemma}
\newtheorem{Corollary}[Theorem]{\sc Corollary}
\newtheorem{Proposition}[Theorem]{\sc Proposition}

\theorembodyfont{\normalfont}
\newtheorem{Remark}[Theorem]{\sc Remark}
\newtheorem{Acknowledgements}[Theorem]{\sc Acknowledgements}

\newtheorem{Definition}[Theorem]{\sc Definition}

\newtheorem{Proof}{\sc Proof}

\makeatletter
\def\im{\mathop{\operator@font im}\nolimits}
\def\red{{\mathop{\operator@font red}\nolimits}}
\def\arg{\mathop{\operator@font arg}\nolimits}

\def\pole{\mathop{\operator@font pole}\nolimits}
\def\tor{{\mathop{\operator@font tor}\nolimits}}
\def\Cst{\mathop{\operator@font Cst}\nolimits}
\def\Newton{\mathop{\operator@font Newton}\nolimits}
\def\Hom{\mathop{\operator@font Hom}\nolimits}
\def\Ext{\mathop{\operator@font Ext}\nolimits}
\def\End{\mathop{\operator@font End}\nolimits}
\def\Ass{\mathop{\operator@font Ass}\nolimits}
\def\Spec{\mathop{\operator@font Spec}\nolimits}
\def\Supp{\mathop{\operator@font Supp}\nolimits}
\def\rk{\mathop{\operator@font rk}\nolimits}
\def\ker{\mathop{\operator@font ker}\nolimits}
\def\dim{\mathop{\operator@font dim}\nolimits}
\def\prof{\mathop{\operator@font depth}\nolimits}
\def\pd{\mathop{\operator@font pd}\nolimits}
\def\coker{\mathop{\operator@font coker}\nolimits}

\def\m{\mathop{\operator@font m}\nolimits}

\def\C{\textbf{\upshape C}}
\def\F{\textbf{\upshape F}}
\def\K{\textbf{\upshape K}}
\def\H{\textbf{\upshape H}}
\def\U{\textbf{\upshape U}}
\def\V{\textbf{\upshape V}}

\newcommand{\htt}{\mathop{\operator@font ht}\nolimits}
\def\Ann{\mathop{\operator@font Ann}\nolimits}
\def\MCM{\mathop{\operator@font MCM}\nolimits}
\def\FID{\mathop{\operator@font FID}\nolimits}
\def\MC{\mathop{\operator@font MC}\nolimits}
\def\CEC{\mathop{\operator@font CEC}\nolimits}
\def\gr{\mathop{\operator@font gr}\nolimits}
\def\Tor{\mathop{\operator@font Tor}\nolimits}
\def\ext{\mathop{\operator@font ext}\nolimits}

\makeatother

\begin{document}

\baselineskip=13pt

\pagestyle{empty}

\ \vspace{1.7in}

\noindent {\LARGE\bf {An Auslander-Buchsbaum identity 
		      for semi\-dualizing modules}}

\bigskip

\noindent  \ Jan R. Strooker, \  Mathematisch Instituut,
Universiteit Utrecht, Postbus 80010, 3508 TA Utrecht, Nederland. {\it E-mail}: {\
strooker@math.uu.nl

\vspace{2.4cm}

\section{Introduction \hfill\break}

Semidualizing modules were introduced by Foxby \cite{F} and Golod \cite{G}
under another name. They can be described as a poor man's projectives,
and of late occasion a spate of papers, going in different directions, e.g.
\cite{Ar}, \cite{Ch}, \cite {HJ}, \cite{FS}, \cite{SW}, \cite{W}, \cite{HW}. There are several reasons why such modules are of interest, but a broad question is where they behave like projectives and where not. However, many studies focus on a more sophisticated concept, semidualizing complexes. Another important feature is to generalize the theory of Gorenstein dimension of Auslander-Bridger \cite{AB} by replacing projective modules with semidualizing ones and then resolving by modules which are totally reflexive with respect to a specific such
module. A concrete understanding in specific cases of semidualizing modules which are not projective, is in its infancy. 

Let us therefore assume we have a fixed such module at our disposal, and derive a few simple consequences. In this note we prove an Auslander-Buchsbaum identity for this module. The original one holds for modules of finite projective dimension \cite[Cor. 7.1.5]{St}, \cite[Th. 1.3.3]{BH}. We show that a theory based on semidualizing modules runs largely parallel, and this is perhaps the main insight on offer here. We do not claim much originality, just want to
show how much can already be achieved in the classical case, finitely generated modules over a local noetherian ring, by elementary means.

\begin{Acknowledgements} A few of the facts and directions pursued here came
to mind during an exchange of e-mails with Diana White, University of Nebraska 
at Lincoln, for which I thank her. As more often, discussions with my Utrecht
colleague Wilberd van der Kallen helped to clarify certain issues.
\end{Acknowledgements}

\bigskip

\section {Semidualizing modules} \hfill\break} 

Throughout this note, $A$ stands for a fixed commutative noetherian local ring
over which we only consider finitely generated modules. This section deals with 
two facts I learned from Diana White. In \cite[Lemma 2.4]{W}, \cite[Th. 2.5]{W} she proves more in greater generality. For our classical case, we provide
easy proofs. 

\begin{Definition} {An $A$-module $C$ is called semidualizing if
	
	(i)	The natural ringhomomorphism $A \rightarrow \End_{A}(C)$
		is an isomorphism;

	(ii)	$\Ext_{A}^{i}(C,C) = 0$ for every $i>0$.}
\end{Definition}

Clearly the ring $A$ is a dualizing module over itself. In what follows, $C$ shall always denote a fixed semidualizing $A$-module. We shall show that 
conditions (i) and (ii) have a strong influence on the additive, covariant
left exact functor $h_{C} = \Hom(C,-)$ where we drop the subscript $A$. But
first we use its contravariant counterpart $h^{C} = \Hom(-,C)$ to prove 

\begin{Lemma} If $C^{n} = U \oplus V$ for $n \geq 0$, then $U = C^{p}$ and
$V = C^{q}$ with $p+q = n$.
\end{Lemma} 
Here and elsewhere, given an arbitrary module $X$, we write $X^{n}$ for the direct sum of $n$ of its copies. 
\begin{Proof} Applying $h^{C}$ to the direct sum, we find $A^{n} = h^{C}(U) 
\oplus h^{C}(V)$. Over the local ring $A$ projectives are free, and we necessarily have a decomposition $A^{n} = A^{p} \oplus A^{q}$ with 
$p+q=n$. Apply $h^{C}$ once more and write $D = h^{C}\circ h^{C}$. Then $D$ is a covariant additive functor which preserves direct sums. The evaluation map from $C^{n}$ to $D(C^{n}) = C^{n}$ is the identity, which proves the lemma. 
\end{Proof}

\begin{Lemma} Any surjection $\pi: C^{n} \rightarrow C^{q}$ splits and $C^{n} =
C^{p} \oplus C^{q}$.
\end{Lemma}

\begin{Proof} Putting $k$ for the residue class field of $A$, we see that 
$k \otimes \pi$ is a surjection of finite dimensional vectorspaces, so
$q \leq n$. The matrix $T$ which describes $\pi$ is a $q \times n$-matrix with entries in $\End(C) = A$. By Nakayama, not every entry in this matrix can be in the maximal ideal of $A$, so there must be a unit somewhere. Permuting the direct sums $C^{n}$ and $C^{q}$ we can get a unit in the upper left hand corner of the matrix, and by the usual elementary transformations obtain a matrix with 1 in the upper left hand corner and zeroes everywhere else in the first row and column. Now $\pi$ splits off the first copy of $C^{n}$ and $C^{q}$. We continue in this way until we obtain a matrix whose first $q$ columns shape a $q \times q$-identity matrix, which furnishes the desired splitting.
\end{Proof}

\bigskip

\section{Further consequences of condition (i)}

In this short section, we concentrate on the functor $h_{C}$ and derive a few  properties which condition (i) confers on it.

A finite exact sequence is called split exact when it is spliced together from
a number of split short exact squences. It is therefore evident that the 
additive functor $h_{C}$ takes such a sequence where the modules are finite 
direct sums of copies of $C$, to a similar split exact sequence where the $C$'s
are replaced by $A$'s. Suppose however that a finite complex of modules $C^{n}$
is taken to a split exact sequence of frees. Then the original complex was split exact. Indeed, in both sequences the maps are defined by rectangular matrices $T$ with entries in $A$. If there is another rectangular matrix $S$ with $ST$ or $TS$ is the appropriate identity matrix, then $S$ defines a splitting equally well in both cases. 

Here is another property of $h_{C}$. For $x \in A$ and $M$ a module, put $K$
for the kernel of the scalar multiplication $x: M \rightarrow M$. Then 
\begin{displaymath}
0 \rightarrow h_{C}(K) \rightarrow h_{C}(M) \stackrel{x} \rightarrow h_{C}(M) 
\end{displaymath}
is exact.

Looking at associated primes and support \cite[Ch. 4-1]{Bo}, one has    	$\Ass h_{C}(M) = \Ass M \cap \Supp C$ for any $M$. Taking $M = C$, one finds $\Ass A = \Ass C$. Thus both these modules have as support $\Spec A$, so have the same dimension. Moreover $\Ass h_{C}(K) = \Ass K$ so that either both these modules are $0$ or neither. This means that $x$ is a nonzerodivisor in $A$ if and only if it is one on $C$. This argument more generally shows that an arbitrary module $M$ and $h_{C}(M)$ have the same nonzerodivisors, as well as the same associated primes, support and dimension. 

It has been pointed out by S. Sather-Wagstaff that these facts and some of the ones to follow, are already collected as ``basic facts" in \cite[p. 68]{G}. We shall just prove and use what we need as we go along, with our single semidualizing module $C$ and not all modules which are totally reflexive with respect to $C$ as in Golod.  

\bigskip

\section{Condition (ii) comes into play}

Let
\begin{displaymath} 
\C = 0 \rightarrow C_{s}\rightarrow C_{s-1} \rightarrow \cdots \rightarrow C_{1} \rightarrow C_{0} \rightarrow Y \rightarrow 0 
\end{displaymath}
be a resolution of the module $Y$ where each $C_{i}$ is the direct sum of a finite number of copies of $C$. We shall prove that the functor $h_{C}$ takes this into a free resolution of $h_{C}(Y)$. 

	If $Y$ has such a resolution, we say that $Y$ has finite $C$-dimension and $C$-dim $Y$ is then the smallest such $s$ with $C_{s} \neq 0$. 
\begin{Lemma} For every $Y$ of finite $C$-dimension, $\Ext^{i}(C,Y) = 0$
for every $i>0$. Moreover the evaluation map on $Y$ with respect to $C$ is the 
identity.
\end{Lemma}
\begin{Proof} Since the property is clear if $C$-dim $Y$ = 0, suppose 
$0 \rightarrow C_{1} \rightarrow C_{0} \rightarrow Y \rightarrow 0$ is
short exact. Applying $h_{C}$ to this gives a long exact sequence in which (ii) entails that $\Ext^{i}(C,Y) = 0$ for every $i>0$. By induction on the $C$-dimension we verify that this is true for every $Y$ of finite $C$-dimension.
The proof of the last statement, that $Y$ is $C$-reflexive, is also by induction and left to the reader.
\end{Proof}
\cite[Cor. 6.5]{HW} contains a sweeping generalization of the Lemma. Working with Auslander and Bass classes, as introduced by Foxby, allows many of the considerations in this note, and more, to be treated in a wider context. 
The paper just cited provides a convenient framework for this. Still, we persist with our down to earth treatment of the classical case.
\begin{Corollary} The exactness of any short exact sequence of modules of 
finite $C$-dimension is preserved by $h_{C}$.  
\end{Corollary}
In fact, this conclusion remains true for any short exact sequence in which 
$\Ext^{1}(C,X) = 0$ for the left hand module $X$. 
\begin{Theorem} Let $\C$ be a resolution of $Y$ as above. Then $h_{C}$
takes this to a free resolution $\F$ of $h_{C}(Y)$. Moreover the projective
dimension $\pd h_{C}(Y) = C\textrm-\dim Y$.
\end{Theorem}
\begin{Proof} Since $\C$ consists of a number of short exact sequences as 
in the previous corollary spliced together, the first statement of the 
theorem follows. Thus $\pd h_{C}(Y) \leq C$-$\dim Y$. According to a classical,
half a century old, result of Eilenberg \cite[Th. 8]{E}                        , the resolution $0 \rightarrow F_{s} \rightarrow F_{s-1} \rightarrow \cdots \rightarrow F_{0} \rightarrow h_{C}(Y) \rightarrow 0$, the $F_{i}$'s finite frees, can be written as a direct sum of a minimal free resolution $\K$ of $h_{C}(Y)$ and a split exact complex $\H$ of frees. Minimal
means that the boundary maps of $\K$ vanish when tensored with the residue 
class field $k$, and then $\pd h_{C}(Y)$ is the length $t$ of $\K$. Since $F_{i} = K_{i} \oplus H_{i}$ for every $i$, the argument in Lemma 2.2 and the fact that the boundary maps in $\C$ and in $\F$ are described by the same matrices with entries in $A$, show that $\C$ = $\U \oplus \V$ where the first complex is taken to $\K$ by $h_{C}$ and the second to $\H$. By a remark in section 3, the complex $\V$ is itself split, and we remove it to find $\U$ as a $C$-resolution of $Y$. This complex has length $t \leq s$. Since $s = C$-$\dim Y$, one must have $s=t$.
\end{Proof}
The first part of this theorem is also a consequence of \cite[Th. 1]{HW}. Diana
White informs me that the second part was proved by her as part of her thesis work and is to appear in a joint preprint with Ryo Takahashi. Consider the 
present proof as a lowbrow approach.  

	To prepare for the next section we make a quick observation. Let 
\begin{displaymath} 
\K = \cdots \rightarrow K_{s} \rightarrow K_{s-1} \rightarrow \cdots \rightarrow 
K_{1} \rightarrow K_{0} \rightarrow 0
\end{displaymath}
be a complex whose homology is concentrated in degree 0. Suppose scalar 
multiplication by $x \in A$ is injective on every chain module $K_{i}$,
and also on $H_{0}(\K)$. Then the factor complex $\K/x\K$ has its homology concentrated in degree 0 where it is $H_{0}(\K)/xH_{0}(\K)$. Indeed, the short exact sequence of complexes gives rise to a long exact sequence of homology modules. The exactness of $\K$ in all positive degrees shows that $\K/x\K$ is 
exact in all degrees greater than 1. At degree 1 exactness follows from 
the injectivity of $x$ on $H_{0}(\K)$.

\bigskip

\section{Setting up for induction}

Let $x \in A$ be a nonzerodivisor, which is automatically a nonzerodivisor
on $C$, see section 3. We want to show that $C/xC$ is a semidualizing module 
over the local ring $A/(x)$.

First notice that multiplying the isomorphism $A \rightarrow \End(C)$ with 
$x$, we get an isomorphism $A/(x) \rightarrow \End(C)/x\End(C)$. Now
applying $h_{C}$ to the short exact sequence $0 \rightarrow C \stackrel{x}
\rightarrow C \rightarrow C/xC \rightarrow 0$, we obtain in view of Corollary
4.2, another short exact sequence which teaches us that $h_{C}(C/xC) =
\End(C)/x\End(C)$ and that every homomorphism $f: C \rightarrow C/xC$ 
factors through $C$. Together this shows that the map from $A/(x)$ to 
$\End_{A/(x)}(C/xC)$ is an isomorphism. 

Next we need to see that $\Ext^{i}_{A/(x)}(C/xC,C/xC) = 0$ for $i >0$.
Take a free $A$-resolution $\K$ of $C$. Then multiplying with $x$ is 
injective on every free chain module $K_{i}$ and $x$ is also injective 
on $H_{0}(\K) = C$. By the remark which ends the previous section, 
we know that $\K/x\K$ provides a free resolution of the $A/(x)$-module
$C/xC$. Now $\Ext^{i}_{A/(x)}(C/xC,C/xC) = H^{i}(\Hom_{A/(x)}(\K/x\K,C/xC)) = 
H^{i}(\Hom_{A}(\K,C/xC)) = \Ext^{i}_{A}(C,C/xC)$. The latter is 0 for all 
$i>0$; this follows from Lemma 4.1 since $C$-dim $C/xC$ = 1. We have proved

\begin{Proposition} Let $C$ be a semidualizing module over the noetherian 
local ring $A$. For every nonzerodivisor in $A$, the module $C/xC$
is semidualizing over $A/(x)$.
\end{Proposition}

Let $Y$ be a module of finite $C$-dimension, and suppose $x$ is a nonzerodivisor  on $Y$ and on $C$.  For the complex $\K$ at the end of the previous section, 
choose a finite $C$-resolution of $Y$. Going through the paces, one sees that
$Y/xY$ has a finite $C/xC$-dimension as an $A/(x)$-module. 

Applying Corollary 4.2 to the short exact $0 \rightarrow Y \stackrel{x}
\rightarrow Y \rightarrow Y/xY \rightarrow 0$ one gets $h_{C}(Y/xY) = h_{C}(Y)/xh_{C}(Y)$, and the left hand term in this identity is $\Hom_{A/(x)}(C/xC,Y/xY)$.

\bigskip

\section{Main result} 

\begin{Theorem} Let $C$ be a semidualizing module over the local noetherian
ring $A$, and let $Y$ be a module of finite $C$-dimension $s$. Then 
$C$-$\dim Y = \prof C - \prof Y$.
\end{Theorem}

\begin{Proof} Let $\C$ be a $C$-resolution of $Y$ of length $s$.
We have seen in Theorem 4.3 that $h_{C}(Y)$ has a minimal free resolution of the same length $s$, say $\F$. Now let $x$ be a nonzerodivisor on $Y$, by section 3 it is also one on $h_{C}(Y)$. According to \cite[section 13.1.9]{St}
and \cite[Th. 9.4.7, Rem. 9.4.8]{BH} it is even one in $A$, because the module $h_{C}(Y)$ has finite projective dimension. Thus $x$ is a nonzerodivisor on $C$. The previous section tells us that over the local ring $A/(x)$, we are in the same position. Namely, $C/xC$ is semidualizing over this ring, and $h_{C}(Y)/xh_{C}(Y) =\Hom_{A/(x)}(C/xC,Y/xY)$. By now it is standard that 
$\C/x\C$ is a $C/xC$-resolution of $Y/xY$ and there is no shorter one.

If we take the minimal resolution $\F$ for $\K$ at the end of section 4,
the free $A/(x)$-resolution $\F/x\F$ of $\Hom_{A/(x)}(C/xC,Y/xY)$
is also minimal, so we recover the well-known fact that 
$\pd_{A/(x)} h_{C}(Y)/xh_{C}(Y) = \pd_{A} h_{C}(Y)$ \cite[Ch. 4-1, Th. E]{Ka}. The left hand term is equal to the $C/xC$-dimension of the $A/(x)$-module $Y/xY$ says Theorem 4.3. 

We continue dividing out by a regular sequence $x_{1}, ..., x_{t}$ on $Y$ 
where $t = \prof Y$, $x_{1} = x$, and write $\mathfrak a$ for the ideal it generates. Over the factor ring $A/\mathfrak a$ we have $\prof Y/\mathfrak a Y = \prof h_{C}(Y)/\mathfrak a h_{C}(Y) = 0$, the latter module still having projective dimension $\pd h_{C}(Y)$. By the Auslander-Buchsbaum formula this integer is $\prof A/\mathfrak a = \prof A - t$, so that $\pd h_{C}(Y) = \prof A - \prof h_{C}(Y)$. Carrying this over to the other side, we find 
$C$-$\dim Y = \prof C - \prof Y$.
\end{Proof}
In semidualizing theory, Auslander-Buchsbaum identities have been stated
before. For instance, working with total reflexivity, as basic fact 8 in \cite{G}, and for semidualizing complexes in \cite[Th. 3.14]{Ch}. From Diana White I understand that she obtained a result which implies the above version, to appear in the joint paper with Takahashi.

We summarize a number of identities which stem from the affinity
between our ring $A$ and its semidualizing module $C$, and the close fit 
between $C$-resolutions of a module $Y$ and the free resolutions of $h_{C}(Y)$
into which $h_{C}$ takes them, in

\begin{Corollary} Let $C$ be a semidualizing module over the noetherian local ring $A$. Then $\Ass C = \Ass A, \Supp C = \Spec A, \dim C = \dim A$, and
both modules have the same nonzerodivisors and the same depth.
For any module $Y$, one has $\Ass Y = \Ass h_{C}(Y), \Supp Y = \Supp h_{C}(Y), 
\dim Y = \dim h_{C}(Y)$, and both modules have the same nonzerodivisors. 
If $C$-$\dim Y$ is finite, it is equal to $\pd h_{C}(Y)$ and $\prof Y = \prof h_{C}(Y)$.
\end{Corollary}

\begin{Remark} As the reader will have noticed, the arguments have all been immediate and standard from the definition of semidualizing onward. Until the
proof of the last theorem. Here we used the fact that a nonzerodivisor on
a module of finite projective dimension is a nonzerodivisor in the ring. 
Or even that a regular sequence on such a module is a regular sequence in the
ring. This was conjectured by M. Auslander \cite[Ch. 2-3]{PS} and essentially proven there in the equal characteristic case, and then in \cite{Ro} for all rings. It belongs to the so-called ``Homological Conjectures" which apparently operate at a deeper level. This one generalizes slightly to the final result. Another generalization was obtained by my then-time student J. Bartijn for equicharacteristic rings \cite[Th. 13.1.9]{St}. 
\end{Remark} 

\begin{Corollary} A regular sequence on a module of finite $C$-dimension
extends to a maximal regular sequence on $C$ which is also regular in the ring. 
\end{Corollary}

\end{document}